\documentclass[a4paper,oneside,12pt]{amsart}
\usepackage{amsmath,amstext,amsthm,a4,amssymb,amscd}
\usepackage[english]{babel}
\usepackage{geometry}
\usepackage{fancyhdr}
\usepackage{fouriernc}
\usepackage{charter}
\usepackage{mathrsfs,dsfont}
\newcommand{\comment}[1]{}
\title{Equidistribution of zeros of holomorphic sections in the non compact setting}
\author{Tien-Cuong Dinh}
\address{UPMC Uni Paris 06, UMR 7586, Institut de
Math{\'e}matiques de Jussieu, 4 place Jussieu, F-75005 Paris,
France}
\email{dinh@math.jussieu.fr }
\thanks{Partially supported by a DFG funded Mercator Professorship}
 \urladdr{www.math.jussieu.fr/{~}dinh}
\author{George Marinescu}
\address{Universit{\"a}t zu K{\"o}ln,  Mathematisches Institut,
    Weyertal 86-90,   50931 K{\"o}ln, Germany\\
    \& Institute of Mathematics `Simion Stoilow', Romanian Academy,
Bucharest, Romania}
\thanks{Partially supported by SFB TR 12}
\email{gmarines@math.uni-koeln.de}
\urladdr{www.mi.uni-koeln.de/{~}gmarines}
\author{Viktoria Schmidt}
\address{Universit{\"a}t zu K{\"o}ln,  Mathematisches Institut,
    Weyertal 86-90,   50931 K{\"o}ln, Germany}
\thanks{Supported by Graduiertenkolleg 1269
Globale Strukturen in Geometrie und Analysis}
\email{vschmidt@math.uni-koeln.de}


\newtheorem{theorem}{Theorem}[section]
\newtheorem{lemma}[theorem]{Lemma}
\newtheorem{proposition}[theorem]{Proposition}
\newtheorem{corollary}[theorem]{Corollary}
\theoremstyle{definition}
\newtheorem{remark}[theorem]{Remark}
\theoremstyle{definition}
\newtheorem{definition}[theorem]{Definition}

\newtheorem{example}[theorem]{Example}
\numberwithin{equation}{section}
\DeclareMathOperator{\Div}{Div}
\DeclareMathOperator{\aut}{Aut}
\DeclareMathOperator{\End}{End}
\DeclareMathOperator{\ric}{Ric}
 \def\cC{\mathscr{C}}
 \def\cO{\mathscr{O}}
  
\newcommand{\calig}[1]{\mathcal{#1}}
\newcommand\mO{\calig{O}}
\newcommand{\field}[1]{\mathbb{#1}}
\newcommand{\Z}{\field{Z}}
\newcommand{\R}{\field{R}}
\newcommand{\C}{\field{C}}
\newcommand{\N}{\field{N}}

\newcommand{\Proj}{\field{P}}
\newcommand{\D}{\field{D}}
\newcommand{\norm}[1]{\lVert#1\rVert}
\newcommand{\FS}{{\rm FS}}
\newcommand{\ddc}{{dd^c}}
\newcommand{\db}{\overline\partial}
\newcommand{\imat}{\sqrt{-1}}
\newcommand{\ov}{\overline}
\newcommand{\wi}{\widetilde}
\newcommand{\boldsym}[1]{\boldsymbol{#1}}
\newcommand\bb{\boldsym{b}}
\begin{document}
\begin{abstract}
We consider tensor powers $L^{N}$ of a positive Hermitian line bundle $(L,h^L)$ over a non-compact complex manifold $X$. In the compact case, B.\ Shiffman and S.\ Zelditch proved that the zeros of random sections become asymptotically uniformly distributed as $N \to \infty$ with respect to the natural measure coming from the curvature of $L$. Under certain boundedness assumptions on the curvature of the canonical line bundle of $X$ and on the Chern form of $L$ we prove a non-compact version of this result. We give various applications, including the limiting distribution of zeros of cusp forms with respect to the principal congruence subgroups of $SL_2(\mathbb Z)$ and to the hyperbolic measure, the higher dimensional case of arithmetic quotients and the case of orthogonal polynomials with weights at infinity.
We also give estimates for the speed of convergence of the currents of integration on the zero-divisors.
\end{abstract}
\maketitle

\tableofcontents

\maketitle

\section{Introduction and related results}
This note is concerned with the asymptotic distribution of zeros of random holomorphic sections in the high tensor powers $L^{N}$ of a positive Hermitian line bundle $(L,h^L)$ over a non-compact complex manifold $X$.
Distribution of zeros of random polynomials is a classical subject, starting with the papers of Bloch-P\'olya, Littlewood-Offord, Hammersley, Kac and Erd\"os-Tur\'an, see e.g. Bleher-Di \cite{BDi:97} and
Shepp-Vanderbei \cite{SheppVanderbei} for a review and complete references.

Shiffman and Zelditch \cite{SZ99} obtained a far-reaching generalization by proving
that the zeros of random sections of powers $L^{N}$ of a positive line bundle $L$ over a projective manifold become asymptotically uniformly distributed as $N \to \infty$ with respect to the natural measure coming from $L$. In a long series of papers these authors further considered the correlations between zeros and their variance (see e.g.\ \cite{BSZ,SZ08}). Berman \cite{Ber10} generalized some of these results  in the context of pseudoconcave domains.

Ideas from complex dynamics were also applied in this field.  Notably, 
a new method to study the distribution of zeros was introduced by Sibony and the first-named author \cite{DS:06} using the formalism of meromorphic transforms. They also gave bounds for the convergence speed in the compact case, improving on the ones in \cite{SZ99}.
Other methods due to Forn\ae ss-Sibony \cite{FS95b} found application
in the study of the distribution of zeros of quantum ergodic eigenfunctions \cite{SZ99}. 

There is an interesting connection between equidistribution of zeros and Quantum Unique Ergodicity related to a conjecture of Rudnik and Sarnak \cite{RuSa:94} about the behaviour of high energy Laplace eigenfunctions on a Riemannian manifold. By replacing Laplace eigenfunctions with modular forms one is lead to study of the equidistribution of zeros of Hecke modular forms. This was done by Rudnick \cite{Rud05}, Holowinsky and Soundararajan \cite{HolSou:10} and generalized by Marshall \cite{Mars11} and Nelson \cite{Nel:10}.

Another area where random polynomials and holomorphic sections play a role is statistical physics. Holomorphic random sections provide a model for quantum chaos and the distribution of their zeros was intensively studied by physicists e.g.\ \cite{BDi:97,BBL,Han:96,NoVo:98,SoTsi04}.

The proof of the equidistribution in \cite{DS:06,SZ99} involves the asymptotic expansion of the Bergman kernel. In \cite{MM07,MM08a} we obtained the asymptotic expansion of the $L^2$-Bergman kernel for positive line bundles over complete Hermitian manifolds under some natural curvature conditions   \eqref{i} (see Remarks \ref{rem1}, \ref{rem-forms}). In this paper we regain the asymptotic equidistribution of random zeros of holomorphic $L^2$-sections under the curvature assumption \eqref{i} and a growth condition \eqref{upb} on the spaces of holomorphic $L^2$-sections of $L^N$.

Let us consider an $n$-dimensional complex Hermitian manifold $(X,J,\Theta)$, where $J$ is the complex structure and $\Theta$ a positive $(1,1)$-form. The manifold $(X,J,\Theta)$ is called K\"ahler if $d\Theta=0$. To $\Theta$ we associate a $J$-invariant Riemannian metric $g^{TX}$ given by $g^{TX}(u,v)=\Theta(u,Jv)$
for all $u,v\in T_xX$, $x\in X$.

We consider further a Hermitian holomorphic line bundle $(L,h^L)$ on $X$. The curvature form of $L$ is denoted by $R^L$. We denote by $L^N:=L^{\otimes N}$ the $N$-th tensor power of $L$. The Hermitian metrics $\Theta$ and $h^L$ provide an $L^2$ Hermitian inner product on the space of sections of $L^N$ and we can introduce the space of holomorphic $L^2$-sections $H^0_{(2)}(X,L^N)$, cf.\ \eqref{lm2.0}.

For a section $s\in H^0_{(2)}(X,L^N)$ we denote by $\Div(s)$ the divisor defined by $s$; then $\Div(s)$ can be written as a locally finite linear combination $\sum c_iV_i$, where $V_i$ are irreducible analytic hypersurfaces and $c_i\in\Z$ are the multiplicities of $s$ along the $V_i$'s.

Recall here the notion of
current: let $\Omega^{p,\,q}_0(X)$ denote the space of smooth compactly supported $(p,q)$-forms on
$X$, and we let $\Omega^{\prime\,p,q}(X) =\Omega_0^{n-p,\,n-q}(X)^{\prime}$ denote the space of
$(p,q)$-currents on $X$; $(T, \varphi) = T(\varphi)$ denotes the pairing
of $T\in\Omega^{\prime\,p,q}(X)$ and $\varphi\in\Omega_0^{n-p,\,n-q}(X)$.

 We denote by $[\Div(s)]$ the current of integration on $\Div(s)$. If $\Div(s)=\sum c_iV_i$ then
 \[
 \big([\Div(s)],\varphi\big):=\sum_ic_i\int_{V_i}\varphi\,,\quad \varphi\in\Omega_0^{n-1,\,n-1}(X)\,,
 \]
 where the integrals are well-defined by a theorem of Lelong (cf.\ \cite[III-2.6]{De:01}, \cite[p.\,32]{GH:78}).

Assume that the spaces $H^0_{(2)}(X,L^N)$ have finite dimension $d_{N}$ for all $N$.
We endow $H^0_{(2)}(X,L^N)$ with the natural $L^{2}$-metric (cf.\ \eqref{0c2}) and this induces a
Fubini-Study metric $\omega_{\FS}$ on the projective space $\Proj H^0_{(2)}(X,L^N)$. The volume form $\omega_{\FS}^{d_{N}-1}$ defines by normalization a probability measure $\sigma_{\FS}$ on $\Proj H^0_{(2)}(X,L^N)$.
We consider the probability space
\begin{equation}\label{probsp}
(\Omega,\mu)=\prod_{N=1}^\infty (\Proj H^0_{(2)}(X,L^N),\sigma_{\FS})\,.
\end{equation}

Note that for two elements $s,s'\in H^0_{(2)}(X,L^N)$ which are in the same equivalence class in $\Proj H^0_{(2)}(X,L^N)$ we have $\Div(s)=\Div(s')$, so $\Div$ is well defined on $\Proj H^0_{(2)}(X,L^N)$.
\begin{definition}
We say that the zero-divisors of generic random sequences $(s_{N})$ with $s_{N}\in \Proj H^0_{(2)}(X,L^N)$ are equidistributed with respect to a Hermitian metric $\omega$ on $X$, if
for $\mu$-almost all sequences $(s_{N}) \in \prod_{N=1}^\infty \Proj H^0_{(2)}(X,L^N)$ we have
\begin{equation}
\frac1N\big[\Div(s_N)\big]\to\omega\,,\quad N\to\infty\,,
\end{equation}
in the sense of currents, that is, for any test $(n-1,n-1)$-form $\varphi\in\Omega_0^{n-1,\,n-1}(X)$
there holds
\begin{equation}
\label{conveq}\lim_{N\to \infty}\Big(\frac{1}{N}\big[\Div(s_N)\big],\varphi\Big)= \int_X \omega\wedge\varphi.
\end{equation}
\end{definition}
Our first result is a generalization to non-compact manifolds of a seminal result of Shiffman-Zelditch \cite[Th.\,1.1]{SZ99}.

Let $K_X$ denote the canonical bundle of $X$, i.e., $K_X=\Lambda^{n,0}T^*X$. If $\Theta$ is a Hermitian metric on $X$ we consider the induced metric on $K_{X}$ with curvature $R^{K_X}$.
If $\Theta$ is K\"ahler, then $\imat R^{K_{X}}=-\ric_\Theta$, where $\ric_\Theta$ is the Ricci curvature of the Riemannian metric associated to $\Theta$ (cf.\ \eqref{ric0}-\eqref{ric1}).

\begin{theorem}\label{conv}
Let $(X,\Theta)$ be an $n$-dimensional complete Hermitian manifold. Let $(L,h^L)$ be a Hermitian holomorphic line bundle over $X$. Assume that there exist constants $\varepsilon>0$, $C>0$ such that
\begin{equation}\label{i}
\imat\,R^L>\varepsilon \Theta,\quad \imat R^{K_X}<C\Theta,\quad |\partial \Theta|_{g^{TX}}<C\,.
\end{equation}
Furthermore, assume that for all $N$ we have $\dim H^0_{(2)}(X,L^N)<\infty$  and 
\begin{equation}\label{upb}
\dim H^0_{(2)}(X,L^N)={\rm O}(N^{n})\,,\quad N\longrightarrow\infty.
\end{equation} 
Then the zero-divisors of generic random sequences $(s_{N}) \in \prod_{N=1}^\infty \Proj H^0_{(2)}(X,L^N)$
are equidistributed with respect to $\frac{\imat}{2\pi}R^L$.
\end{theorem}
\begin{remark}\label{rem1}
Note that \eqref{i} are the natural hypotheses to apply H\"ormander's $L^2$ $\overline\partial$-method (with singular weights) in order to produce $L^2$ holomorphic sections of $L^N$ which separate points and give local coordinates. 

In this paper the existence of sections appears in the guise of the asymptotic expansion of the Bergman kernel of $L^N$ (cf.\ \cite{MM08a}, \cite[Th.\,6.1.1]{MM07}, see Theorem \ref{bke}); this follows from the spectral gap of the Kodaira-Laplacian on $L^N$, which is in turn a consequence of \eqref{i} via the Bochner-Kodaira-Nakano formula

Assumption \eqref{upb} says that the growth of the dimension of the spaces of holomorphic sections is the same as on a compact or more generally a pseudoconcave manifold (see Sections 4-7).

More details about the hypotheses \eqref{i} and other variants (for $n$-forms or for $L=K_X$) are discussed in Remarks \ref{rem-forms}, \ref{conv-var}. 
For the moment let us just observe that if $(X,\Theta)$ is K\"ahler, then $\partial \Theta=0$, so the condition on $\partial\Theta$ in \eqref{i}
is trivially satisfied.  
\end{remark}
\begin{remark}
Note that on a compact manifold equipped with a positive line bundle one can always realize condition \eqref{i}  by choosing $\Theta=\frac{\imat}{2\pi}R^L$. This is not always possible in general; if $X$ is non-compact, the metric associated to $\Theta=\frac{\imat}{2\pi}R^L$ might be non-complete.
The existence of a complete metric $\Theta$ with $\imat R^L>\varepsilon \Theta$ is equivalent to saying that $\imat R^L$ defines a complete K\"ahler metric.

Another interesting fact is that by changing $\Theta$ we change the $L^{2}$-product \eqref{0c2}
on the spaces of sections of $L^{N}$. However, as long as \eqref{i} holds, the $L^{2}$-holomorphic sections are equidistributed with respect to $\tfrac{\imat}{2\pi}R^L$.
\end{remark}
Theorem \ref{conv} is a consequence of the following equidistribution result on relatively compact open sets. We also address here the problem of the convergence speed of the currents $[\frac{1}{N}\Div (s_N)]$ towards $\frac{\imat}{2\pi}R^{L}$. For a compact manifold $X$ the estimates of the convergence speed from Theorem \ref{conv-speed1} were obtained in \cite{DS:06} and we will adapt the method of \cite{DS:06} in the present context.

\begin{theorem}\label{conv-speed1}
Let $(X,\Theta)$ and $(L,h^L)$ be as in Theorem \ref{conv}.
Then for any relatively compact open subset $U$ of $X$ there exist a constant $c=c(U)>0$ and an integer $N(U)$ with the following property. For any real sequence $(\lambda_N)$ with $\lim_{N\to\infty} (\lambda_N/\log N)=\infty$ and for any $N\geqslant N(U)$ there exists a set $E_{N}\subset\mathbb{P}H^0_{(2)}(X,L^{N})$ such that:

\begin{enumerate}
\item[(a)] $\sigma_\FS(E_{N})\leqslant c\,N^{2n}\,e^{-\lambda_{N}/c}$,
\item[(b)] For any $s\in\mathbb{P}H^0_{(2)}(X,L^{N})\setminus E_{N}$ we have the following estimate
\begin{equation}\label{conv-est1}
\Big|\Big(\frac1N\big[\Div (s)\big]-\frac{\imat}{2\pi}R^{L},\varphi\Big)\Big|\leqslant \frac{\lambda_{N}}{N}\|\varphi\|_{\cC^{2}}\,,\quad \varphi\in\Omega_0^{n-1,n-1}(U)\,.
\end{equation}
\end{enumerate}
In particular, for a generic sequence $(s_{N})\in \prod_{N=1}^{\infty}\mathbb{P}H^0_{(2)}(X,L^{N})$ the estimate \eqref{conv-est1} holds for $s=s_{N}$ and $N$ large enough.
\end{theorem}
By choosing $(\lambda_{N})$ such that $\lim_{N\to\infty}(\lambda_{N}/N)=0$ we obtain
that zero-divisors of generic random sequences $(s_{N}) \in \prod_{N=1}^\infty \Proj H^0_{(2)}(X,L^N)$
are equidistributed with respect to $\frac{\imat}{2\pi}R^L$ on $U$. From this observation, Theorem \ref{conv} follows immediately.

The sequence $(E_N)_{N\geqslant N(U)}$ in Theorem \ref{conv-speed1} depends on the choice of the sequence $(\lambda_N)$. Thus, we can use different choices for $(\lambda_{N})$ depending of our purpose: if we want a better estimate on $E_N$, the speed of equidistribution will be worse; if we want a better convergence speed, we have to allow larger sets $E_N$. For example, we can take $\lambda_{N}=(\log N)^{1+\beta}$, $0<\beta\ll1$. Then $(\lambda_{N}/N)$ in \eqref{conv-est1} converges very fast to $0$. 
Note that in \cite[p.\,671]{SZ99} it is observed that in the compact case the convergence speed is bounded above by $N^{\varepsilon-\frac12}$, for any $\varepsilon>0$.

Let $X$ be an $n$-dimensional, $n\geqslant2$, irreducible  quotient of a bounded symmetric domain
$D$ by a torsion-free arithmetic group $\Gamma\subset\aut(D)$. We call such manifolds arithmetic quotients. The Bergman metric $\omega^{\mathcal B}_{D}$ on $D$ (cf. \eqref{bergmet}) descends to a Hermitian metric on $X$, called Bergman metric on $X$ and denoted by $\omega^{\mathcal B}_{X}$. This metric induces a volume form on $X$ and a Hermitian metric on $K_{X}$ and thus an
$L^{2}$ Hermitian inner product on the sections of $K_{X}^{N}$.
\begin{corollary}\label{equi-arit}
Let $X$ be an $n$-dimensional, $n\geqslant2$, arithmetic quotient.
Let $H^{0}_{(2)}(X,K_{X}^{N})$ be the space of holomorphic sections of $K_{X}^{N}$ which are square-integrable with respect to the $L^{2}$ Hermitian inner product induced by the Bergman metric.
Then the zero-divisors of generic random sequences $(s_{N}) \in \prod_{N=1}^\infty \Proj H^0_{(2)}(X,K_{X}^{N})$
are equidistributed with respect to the Bergman metric $\frac{1}{2\pi}\omega^{\mathcal B}_{X}$ on $X$. Moreover, we have an estimate of the convergence speed on compact sets as in Theorem \ref{conv-speed1}.
\end{corollary}
Theorem \ref{conv} has the following application to the equidistribution of zeros of modular forms.
\begin{corollary}\label{equi-modular}
Let $\Gamma\subset SL_{2}(\mathbb Z)$ be a subgroup of finite index that acts freely on the hyperbolic plane $\mathbb H.$ Consider the spaces of cusp forms $\mathcal S_{2N}(\Gamma)$ as Gaussian probability spaces with the measure induced by the Petersson inner product. Then for almost all random sequences $(f_{2N})$ in the product space $\prod_{N=1}^\infty\mathcal S_{2N}(\Gamma)$ the associated sequence of zeros becomes asymptotically uniformly distributed, i.e.\ for piecewise smooth open sets $U$ contained in one fundamental domain we have
\[ \lim_{N\to \infty} \frac{1}{N}\sharp \Big\{z\in U; f_{2N}(z)=0\Big\}=\frac{1}{2\pi}\operatorname{Vol}\,(U),\] where $\operatorname{Vol}$ denotes the hyperbolic volume.
Moreover, we have an estimate of the convergence speed on compact sets as in Theorem \ref{conv-speed1}.
\end{corollary}
This result concerns \emph{typical} sequences of cusp forms. If one considers Hecke modular forms, by a result of Rudnick \cite{Rud05} and its generalization by Marshall \cite{Mars11}, the zeros of \emph{all} sequences of Hecke modular forms are equidistributed. The method used in \cite{Mars11,Rud05} follows the seminal paper of Nonnenmacher-Voros \cite{NoVo:98} and consists
in showing the equidistribution of masses of Hecke forms. Rudnick \cite{Rud05} invoked for this ¥purpose the Generalized Riemann Hypothesis and this hypothesis was later removed by Holowinsky and Soundararajan \cite{HolSou:10}. Marshall \cite{Mars11} extended their methods to the higher-dimensional setting. Nelson \cite{Nel:10} removed later some of the hypotheses in \cite{Mars11}. Cogdell-Luo \cite{CL08} proved the mass-equidistribution for Hecke eigenforms
on the Siegel modular varieties.

The lay-out of this paper is as follows. In Section \ref{sec:bergker}, we collect the necessary ingredients about the asymptotic expansion of Bergman kernel. In Section \ref{sec:equi_speed}, we prove the main results, Theorems \ref{conv} and \ref{conv-speed1}, about the equidistribution on compact sets together with the estimate of the convergence speed. In the next sections we give applications of our main result in several geometric contexts and prove equidistribution:  for sections of the pluricanonical bundles over pseudoconcave manifolds and arithmetic quotients in dimension greater than two (Section \ref{sec:psconc}), for modular forms over Riemann surfaces (Section \ref{sec:modular}), for sections of positive line bundles over quasi-projective manifolds (Section \ref{sec:quasiproj}) and finally for polynomials over $\C$ endowed with the Poincar\'e metric at infinity (Section \ref{sec:poly}).

\section{Background on the Bergman kernel}\label{sec:bergker}
Let $(X,J,\Theta)$ be a complex Hermitian manifold  of dimension $n$,
where $J$ is the complex structure and $\Theta$ is the $(1,1)$-form associated to a Riemannian metric $g^{TX}$ compatible with $J$, i.e.\
\begin{equation}\label{riem-met}
\Theta(u,v)=  g^{TX}(Ju,v)\,,\: g^{TX}(Ju,Jv)=  g^{TX}(u,v)\,,\quad
\text{for all $u,v\in T_xX$, $x\in X$}\,.
\end{equation}
The volume form of the metric $g^{TX}$ is given by $dv_X=\Theta^n/n!$.
The Hermitian manifold $(X,J,\Theta)$ is called complete if $g^{TX}$ is a complete Riemannian metric.

Let $(E,h^E)$ be a Hermitian holomorphic line bundle. For $v,w\in E_x$, $x\in X$, we denote by
$\langle v,w\rangle_{E}$ the inner product given by $h^E$ and by $|v|_{E}=\langle v,v\rangle_{E}^{1/2}$ the induced norm.
The $L^2$--Hermitian product on the space $\cC^\infty_0 (X,E)$
of compactly supported smooth sections of $E$, is given by
\begin{equation}\label{0c2}
(s_1,s_2) =\int_X\big\langle s_1(x),
s_2(x)\big\rangle_{E}\,dv_X(x)\,.
\end{equation}
We denote by $L^2(X,E)$ the completion of $\cC^\infty_0 (X,E)$ with respect to \eqref{0c2}.
Consider further the space of holomorphic $L^2$-sections of $E$:
\begin{equation}\label{lm2.0}
H^{0}_{(2)}(X,E):=\Big\{s\in L^{2}(X,E) : \overline{\partial}^{E}s=0\Big\}\,.
\end{equation}
Here the condition $\overline{\partial}^{E}s=0$ is taken in the sense of distributions.
Elementary continuity properties of differential operators on distributions show that $H^{0}_{(2)}(X,E)$ is a closed subspace of $L^{2}(X,E)$. Moreover, the hypoellipticity of $\overline{\partial}^{E}$ implies that elements of $H^{0}_{(2)}(X,E)$ are smooth and indeed holomorphic.

For a Hermitian holomorphic line bundle $(E,h^E)$ we denote by $R^E$ its curvature, which is a $(1,1)$-form on $X$.
Given the Riemannian $g^{TX}$ metric on $X$ associated to $\Theta$ as in \eqref{riem-met},
we can identify
$R^E$ to a Hermitian matrix
$\dot{R}^E\in \End(T^{(1,0)}X)$ such that for $u,v\in T^{(1,0)}_xX$,
\begin{equation*}\label{herm3.8}
R^E(u,\overline{v})
= \langle \dot{R}^E (u), v\rangle.
\end{equation*}
There exists an orthonormal basis of $(T^{(1,0)}_xX,g^{TX})$ such that $\dot{R}^E=\operatorname{diag}(\alpha_{1},\ldots,\alpha_{n})$. The real numbers $\alpha_{1},\ldots,\alpha_{n}$ are called the eigenvalues of $R^{E}$ with respect to $\Theta$ at $x\in X$. We denote
\begin{equation*}\label{herm3.81}
\det(\dot{R}^{E})=\prod_{k=1}^{n}\alpha_{k}\,.
\end{equation*}
For a holomorphic section $s$ in a holomorphic line bundle, let $\Div(s)$ be the zero divisor of $s$ and we denote
by the same symbol $\Div(s)$ the current of integration on $\Div(s)$.
We denote by $K_X$ the canonical line bundle over $X$.

Let $(L,h^L)\to X$ be a Hermitian holomorphic vector bundle.
As usual, we will write $L^N$ instead of $L^{\otimes N}$ for the $N-$th tensor power of $L$ and omit the $\otimes-\text{symbol}$ in all similar expressions. On $L^N$ we consider the induced Hermitian metric $h^{L^N}=(h^{L})^{\otimes N}$. With respect to this metrics we define as in \eqref{0c2} an $L^2$-Hermitian product on the spaces $\cC^{\infty}_0(X,L^N)$.
Denote by $H^0_{(2)}(X,L^{N})$ the subspace of holomorphic $L^2$-sections as in \eqref{lm2.0}. The Bergman kernel $P_N(x,y)$ is the Schwartz kernel of the orthogonal projection
\[P_N:L^2(X,L^N)\to H^0_{(2)}(X,L^{N})\,.
\]
If $H^0_{(2)}(X,L^{N})=0$, we have of course $P_N(x,x)=0$ for all $x\in X$. If $H^0_{(2)}(X,L^{N})\neq0$, consider an orthonormal basis $(S_j^N)_{j=1}^{d_N}$ of $H^0_{(2)}(X,L^{N})$ (where $1\leqslant d_N\leqslant\infty$). Then
\begin{equation}\label{bk}
P_N(x,x)= \sum_{j=1}^{d_N} |S_j^N(x)|_{L^N}^2\quad \text{in $\cC_{loc}^\infty(X)$}\,.
\end{equation}
The following proposition is a special case of \cite[Th.\,6.1.1]{MM07}. The proof is based on the existence of the spectral gap \eqref{spec-gap} and the analytic localization techniques of Bismut-Lebeau \cite{BL91}. The case when $X$ is compact is due to Catlin \cite{Catlin} and Zelditch \cite{Zelditch98}.
\begin{theorem}\label{bke}
Let $(X,\Theta)$ be an $n$-dimensional complete Hermitian manifold. Let $(L,h^L)$ be a Hermitian holomorphic line bundle over $X$. Assume that there exist constants $\varepsilon>0$, $C>0$ such that \eqref{i} is fulfilled.
Then there exist coefficients $\bb_r \in \cC^{\infty}(X)$, for $r\in\mathbb N$, such that the following asymptotic expansion
\begin{equation}\label{bergexp0}
P_N(x,x) =\sum_{r=0}^{\infty}\bb_r(x)N^{n-r}
\end{equation}
holds in any $\cC^\ell$-topology on compact sets of $X$. Moreover, $b_0= \det \left(\frac{\dot{R}^L}{2\pi}\right)$.
\end{theorem}
To be more precise, the asymptotic expansion \eqref{bergexp0} means that
for any compact set $K\subset X$ and any $k,\ell \in \mathbb N$ there exists a constant $C_{k,\ell,K}>0$, such that for any $N \in \mathbb N$,
\begin{equation}\label{bergexp}
\left|P_N(x,x) - \sum_{r=0}^{k}\bb_r(x)N^{n-r}\right|_{\cC^\ell(K)}\leq C_{k,\ell,K}N^{n-k-1}\,.
\end{equation}
Let
\begin{equation}\label{baseloc}
Bs_N:=\{x\in X: s(x)=0\;\text{ for all } s\in H_{(2)}^0(X,L^N)\}
\end{equation}
be the base locus of $H_{(2)}^0(X,L^N)$, which is an analytic set.
Assume now as in Theorem \ref{conv} that for all $N$ we have
\[
d_{N}:=\dim\,H^0_{(2)}(X,L^N)<\infty\,.
\]
The Kodaira map is the holomorphic map
\begin{equation}\label{kodmap}
\begin{split}
&\Phi_N:\;X\setminus Bs_{N}\to \mathbb PH^0_{(2)}(X,L^N)^*\,,\\
&x \longmapsto \big\{s \in H^0_{(2)}(X,L^N):\,s(x)=0\big\}\,.
\end{split}
\end{equation}
In this definition we identify the projective space $\mathbb PH^0_{(2)}(X,L^N)^*$ of lines in
$H^0_{(2)}(X,L^N)^*$ to the Grassmannian manifold of hyperplanes in $H^0_{(2)}(X,L^N)$.
To be precise, for a section $s\in H^0_{(2)}(X,L^N)$ and a local holomorphic frame $e_L:U\to L$ of $L$, there exists a holomorphic function $f\in\mathcal{O}(U)$ such that
$s=fe_L^{\otimes N}$. We denote $f$ by $s/e_L^{\otimes N}$.
To $x \in X$ and a choice of local holomorphic frame $e_L$ of $L$, we assign the element $s\mapsto (s/e_L^{\otimes N})(x)$ in $H^0_{(2)}(X,L^N)^*$. If $x \notin Bs_N$, this defines a line in $H^0_{(2)}(X,L^N)^*$ which is by definition $\Phi_N(x)$. By a choice of basis $\{S^N_i\}$ of $H^0_{(2)}(X,L^N)$ it is easy to see that $\Phi_N$ is holomorphic and has the coordinate representation \[x \mapsto [(S^N_0/e_L^{\otimes N})(x),\dots,(S^N_{d_N}/e_L^{\otimes N})(x)].\]\par

As a consequence of the asymptotic expansion \eqref{bergexp0} of the Bergman kernel we obtain:
\begin{corollary}\label{kod-emb}
Under the assumption of Theorem \ref{bke} we have:
\\[2pt]
(i) The ring $\oplus_N H^0_{(2)}(X,L^N)$ separates points and gives local coordinates at each point of $X$.
\\[2pt]
(ii) Let $K\subset X$ be compact. There exists an integer $N(K)$ such that for every $N>N(K)$ we have
$Bs_{N}\cap K=\emptyset$ and $\Phi_{N}$ is an embedding in a neighbourhood of $K$.
\\[2pt]
(iii) $\liminf_{N\to\infty}N^{-n}d_N>0$. Thus $d_N\sim N^n$ for $N\to\infty$.
\end{corollary}
\begin{proof}
Item (i) follows by applying the analytic proof of the Kodaira embedding theorem by the Bergman kernel expansion (see e.g.\ \cite[\S 8.3.5]{MM07}, where symplectic manifolds are considered; the arguments therein extend to the non-compact case due to \cite[Th.\,6.1.1]{MM07}). Item (ii) follows immediately from (i).

The first estimate in (iii) follows from Fatou's lemma, applied on $X$ with
the measure $\Theta^n/n!$
to the sequence $N^{-n}P_N(x,x)$ which converges pointwise to $\bb_0$ on $X$. Note that 
\[
\int_XP_N(x,x)\frac{\Theta^n}{n!}=d_N
\] 
by \eqref{bk}. Hence
\begin{equation} \label{ell6}
\liminf_{N\longrightarrow\infty}N^{-n}d_N\geqslant
\frac{1}{n!}\int_X\Big(\tfrac{\sqrt{-1}}{2\pi}R^L\Big)^n >0.
\end{equation}
The estimate $d_N\sim N^n$ for $N\to\infty$ follows from \eqref{ell6} and \eqref{upb}.
\end{proof}

Let us note a case when the Kodaira map \eqref{kodmap} gives a global embedding.
\begin{proposition}
Under the assumption of Theorem \ref{bke} suppose further  that the base manifold $X$ is $1$-concave (cf.\ Definition \ref{def-ag}).
Then there exists $N_0$ such that for all $N\geqslant N_0$ we have $Bs_N=\emptyset$ and
the Kodaira map \eqref{kodmap} is an embedding of $X$.
\end{proposition}
\begin{proof}
By Corollary \ref{kod-emb} (i) the ring $\oplus_N H^0_{(2)}(X,L^N)$ separates points and gives local coordinates at each point of $X$. Applying the proof of the Andreotti-Tomassini embedding theorem \cite{AnTo:70} to the ring $\oplus_N H^0_{(2)}(X,L^N)$ we obtained the desired conclusion.
\end{proof}

Let $V$ be a finite dimensional Hermitian vector space and let $V^*$ be its dual.
Let $\mO(-1)$ be the universal (tautological) line bundle on $\Proj(V^*)$.
Let us denote by $\mO(1)=\mO(-1)^{*}$ the hyperplane line bundle over the projective space $\mathbb P(V^{*})$.

A Hermitian metric $h^V$ on $V$ induces naturally
a Hermitian metric $h^{V^*}$ on $V^*$, thus it induces a Hermitian metric
$h^{\mO(-1)}$ on $\mO(-1)$, as a sub-bundle of
the trivial bundle $V^*$ on $\Proj(V^*)$.
Let $h^{\mO(1)}$ be the Hermitian metric on $\mO(1)$ induced by $h^{\mO(-1)}$.

For any $v\in V$, the linear map $V^*\ni f\to (f,v)\in\C$
defines naturally a holomorphic section $\sigma_v$ of $\mO(1)$ on
$\mathbb P(V^*)$.
By the definition, for $f\in V^*\smallsetminus \{0\}$,
at $[f]\in \Proj(V^*)$, we have
\begin{align}\label{pb1.2}
|\sigma_v([f])|^2_{h^{\cO(1)}} =|(f,v)|^2/|f|^2_{h^{V^*}}.
\end{align}
For $N\geqslant N(K)$, $\Phi_N:K\longrightarrow\mathbb PH^0(X,L^N)^*$
is holomorphic and the map
\begin{equation} \label{pb1.10}
\begin{split}
&\Psi_N: \Phi_N^*\mO(1)\to L^N,\\
&\Psi_N ((\Phi_N^* \sigma_s)(x))=s(x),
\quad \text{for any $s\in H^0(X, L^N)$}
\end{split}
\end{equation}
defines a canonical isomorphism from $\Phi_N^*\mO(1)$ to $L^N$ on $X$,
and under this isomorphism, we have
\begin{equation}
\label{distortion}
h^{\Phi_N^*\mO(1)}(x)= P_N(x,x)^{-1}h^{L^N}(x)
\end{equation}
on $K$ (see e.g.\ \cite[(5.1.15)]{MM07}). Here $h^{\Phi_N^*\mO(1)}$ is the metric on $\Phi_N^*\mO(1)$ induced by
the Fubini-Study metric $h^{\mO(1)}$ on $\mO(1)\to\mathbb PH^0(X,L^N)^*$.

The Fubini-Study form $\omega_{\FS}$
is the K\"ahler form
on $\Proj(V^*)$, and is defined as follows: for any $0\neq v\in V$ we set
\begin{align}\label{pb1.3}
\omega_{\FS}=\frac{\sqrt{-1}}{2\pi} R^{ \mO(1)}
= \frac{\sqrt{-1}}{2\pi} \overline{\partial}\partial \log |\sigma_v|^2_{h^{\mO(1)}}
\quad \text{on $\big\{x\in \Proj(V^*), \sigma_v(x)\neq 0\big\}$}\,.
\end{align}

Identity \eqref{distortion} immediately implies
\begin{equation}\label{eq1}
\begin{split}
\Phi_N^*\omega_{\FS} &= \frac{\imat}{2\pi}R^{L^N}+\frac{\imat}{2\pi}\partial \overline{\partial}\log P_N(x,x)\\
 &=N\frac{\imat}{2\pi}R^{L} + \frac{\imat}{2\pi}\partial\overline{\partial}\log P_N(x,x).
\end{split}
\end{equation}
We obtain as a consequence the analogue of the Tian-Ruan convergence of the Fubini-Study metric on compact subsets of $X$.
\begin{corollary}[{\cite[Cor.\,6.1.2]{MM07}}]
\label{Tian} In the conditions of Theorem \ref{bke}
for any $\ell\in\N$ there exists $C_{\ell,K}$ such that
\begin{equation}\label{Tianest}
 \left| \frac{1}{N} \Phi_N^*\omega_{\FS} - \frac{\imat}{2\pi}R^{L} \right|_{\cC^{\ell}(K)} \leqslant \frac{C_{\ell,K}}{N}\,\cdot
\end{equation}
Thus the induced Fubini-Study form $\frac{1}{N} \Phi_N^*\omega_{\FS} $ converges to $\omega$ in the $\cC^{\infty}_{loc}$ topology as $N\to \infty$.
\end{corollary}
\begin{proof}
By \eqref{eq1},
\[ \frac{1}{N} \Phi_N^*\omega_{\FS} - \frac{\imat}{2\pi}R^{L} = \frac{\imat}{2\pi N}\partial\overline{\partial}\log P_N(x,x) \]
and by \eqref{bergexp}
\[  \left| \partial\overline{\partial}\log P_N(x,x)\right|_{\cC^\ell(K)}={\rm O}(1)\,,\quad N\to\infty \,,\]
where ${\rm O}(1)$ is the Landau symbol.
\end{proof}
Tian obtained the estimate \eqref{Tianest} for $L=K_X$, with the bound ${\rm O}(1/\sqrt{N})$ and for $\ell=2$, for compact manifolds \cite[Th.\,A]{Tia:90} and for
complete K\"ahler manifolds $(X,\omega)$ such that there exists a constant $k>0$ with $\ric_\omega\leqslant -k\,\omega$, in \cite[Th.\,4.1]{Tia:90}. Ruan \cite{Ruan:98} proved the $\cC^{\infty}$-convergence and  improved the bound to ${\rm O}(1/N)$.
\begin{remark}\label{rem-tianest}
Assume that in Theorem \ref{conv} we have
$\Theta=\frac{\imat}{2\pi}R^L$. Then $\bb_0(x)=1$.
Thus from \eqref{bergexp0} and \eqref{eq1}, we can improve the estimate \eqref{Tianest} by replacing the left-hand side by $C_{\ell,K}/N^{2}$ (cf.\ \cite[Rem.\,5.1.5]{MM07}).
\end{remark}
\begin{remark}\label{rem-forms}
We wish to explain here in more detail the origin of the hypotheses \eqref{i} from Theorems \ref{conv} and \ref{bke}. 
\\[2pt]
(i) Let us denote by $\mathcal{T}=[i(\Theta),\partial\Theta]$ the Hermitian torsion of the Poincar\'e metric $\Theta$. Set $\widetilde{L^N}=L^N\otimes K^*_X$. There exists a natural isometry
\begin{equation*}
\begin{split}
&\Psi=\thicksim\,:\Lambda^{0,q}(T^*X)\otimes L^N
\longrightarrow\Lambda^{n,q}(T^*X)\otimes \widetilde{L^N},
\\&
\Psi \, s=\widetilde s=(w^1\wedge\ldots\wedge w^n\wedge s)\otimes
(w_1\wedge\ldots\wedge w_n),
\end{split}
\end{equation*}
where $\{w_j\}^n_{j=1}$ a local orthonormal frame of
$T^{(1,0)}X$. 
The Bochner-Kodaira-Nakano formula \cite[Cor.\,1.4.17]{MM07} shows that for any $s\in \Omega^{0,1}_0(X,L^N)$ we have
\begin{equation} \label{herm20,111}
\begin{split}
\frac{3}{2}\left(\norm{\db s}^2+\norm{\db^{\,*}s}^2\right)
\geqslant &\,
\left\langle R^{L^N\otimes K^*_X}(w_j,\ov{w}_k)\ov{w}^k\wedge i_{\ov{w}_j}s,
s\right\rangle\\
&-\frac{1}{2}\big(\norm{\mathcal{T}^*\wi{s}}^2
+\norm{\ov{\mathcal{T}}\wi{s}}^2 
+\norm{\ov{\mathcal{T}}^*\wi{s}}^2\big).
\end{split}
\end{equation}
By \eqref{i}, $\imat R^{L^N\otimes K^*_X}=\imat (NR^{L}+R^{K^*_X})\geqslant (N\varepsilon-C)\Theta$ hence 
\[\left\langle R^{L^N\otimes K^*_X}(w_j,\ov{w}_k)\ov{w}^k\wedge i_{\ov{w}_j}s,
s\right\rangle\geqslant (N\varepsilon-C)\norm{s}^2\,,\quad s\in \Omega^{0,1}_0(X,L^N)\,.
\] 
Moreover, the condition $|\partial \Theta|_{g^{TX}}<C$ implies that the torsion operators $\mathcal{T}^*$, $\ov{\mathcal{T}}$, $\ov{\mathcal{T}}^*$ are pointwise bounded. Thus \eqref{herm20,111} shows that there exists a constant $C_1>0$ such that 
\[
\norm{\db s}^2+\norm{\db^{\,*}s}^2\geqslant C_1N\norm{s}^2\,,\quad s\in \Omega^{0,1}_0(X,L^N)\,.
\]
This implies (see \cite[p.\,273]{MM07}) that the spectrum of the Kodaira-Laplace operator $\Box_N=\db^{\,*}\db$ on $L^2(X,L^N)$ has a spectral gap, i.\,e., 
\begin{equation}\label{spec-gap}
\operatorname{Spec}\Box_N\subset\big\{0\big\}\cup \big[C_1N,\infty\big)\,.
\end{equation}
It is the general principle in \cite{MM07} that the spectral gap property implies the asymptotic expansion of the Bergman kernel, see e.\,g.\ \cite[Th.\,6.1.1]{MM07}. 
\\[2pt]
(ii) If we work with holomorphic $n$-forms with values in $L^N$, i.\,e.\ sections of $L^N\otimes K_X$, we need less hypotheses than \eqref{i} in order to obtain the spectral gap. By \cite[Th.\,1.4.14]{MM07}
we have for any $s\in \Omega^{n,1}_0(X,L^N)$
\begin{equation} \label{herm20,1}
\begin{split}
\frac{3}{2}\left(\norm{\db s}^2+\norm{\db^{\,*}s}^2\right)\geqslant\big\langle\big[\sqrt{-1}R^{L^N},
i(\Theta)\big]s,s\big\rangle
-\frac{1}{2}(\norm{\mathcal{T}^*s}^2
+\norm{\overline{\mathcal{T}}s}^2
+\norm{\overline{\mathcal{T}}^*s}^2)\,.
\end{split}
\end{equation}
From \eqref{herm20,1} we deduce as above that if
\begin{equation}\label{i_bis}
\imat\,R^L>\varepsilon \Theta,\quad |\partial \Theta|_{g^{TX}}<C\,,
\end{equation}
then the Kodaira-Laplace operator $\Box_N$ on $L^2(X,L^N\otimes K_X)$ has a spectral gap
and the Bergman kernel of the projection $P_N:L^2(X,L^N\otimes K_X)\to H^0_{(2)}(X,L^N\otimes K_X)$
has an asymptotic expansion as in Theorem \ref{bke}.
\\[2pt]
(iii) Assume that $L=K_X$, $h^{K_X}$ is induced by $\Theta$ and \eqref{i_bis} holds. We can view  $(0,q)$-forms with values in $K_X^N$ as $(n,q)$-forms with values in $K_X^{N-1}$. By (ii), the Bergman kernel of the projection $P_N:L^2(X,K^N_X)\to H^0_{(2)}(X,K^N_X)$
has an asymptotic expansion as in Theorem \ref{bke}. This follows of course also by applying \eqref{herm20,111}.
\end{remark}

\section{Equidistribution on compact sets and speed of convergence}\label{sec:equi_speed}
Let $(X,\omega)$ be a Hermitian complex manifold of dimension $n$. Let $U$ be a relatively
compact open subset of $X$. If $S$ is a
real current of bidegree $(p,p)$ and of order 0 on $X$ and $\alpha \geqslant 0$, define
$$\|S\|_{U,-\alpha}:=\sup_{\varphi}|(S,\varphi)|$$
where the supremum is taken over all smooth real $(k-p,k-p)$-forms $\varphi$ with compact support in $U$
such that  $\|\varphi\|_{\cC^\alpha}\leqslant 1$.
When $\alpha=0$, we obtain the mass of $S$ on $U$ that is denoted by $\|S\|_U$.

It is clear that if $\beta\geq\alpha$, then
$$\|S\|_{U,-\alpha}\geqslant \|S\|_{U,-\beta}\,.$$
If $W$ is an open set such that $U\Subset W\Subset X$, by theory of interpolation between Banach spaces \cite{Triebel95}, we have
$$\|S\|_{U,-\alpha}\leqslant c\|S\|_W^{1-\alpha/\beta}\|S\|_{W,-\beta}^{\alpha/\beta}\,,$$
where $c>0$ is a constant independent of $S$, see \cite{DS09}.

For an arbitrary complex vector space $V$ we denote by $\mathbb{P}(V)$ the projective space of $1$-dimensional subspaces of $V$. Fix now a vector space $V$ of complex dimension $d+1$.
Recall that there is a canonical identification of $\mathbb{P}(V^{*})$ with the Grassmannian $G_{d}(V)$ of hyperplanes in $V$, given by $\mathbb{P}(V^{*})\ni[\xi]\mapsto H_{\xi}:=\ker\xi\in G_{d}(V) $, for $\xi\in V^{*}\setminus\{0\}$.

Once we fix a Hermitian product on $V$, we endow the various projective spaces with the Fubini-Study metric $\omega_{\FS}$, normalized such that the induced measure $\sigma_\FS:=\omega_\FS^{d}$ is a probability measure.

\begin{theorem} \label{th_equi}
Let $(X,\omega)$ be a Hermitian complex manifold of dimension $n$ and let $U$ be a relatively
compact open subset of $X$. Let $V$ be a Hermitian complex vector space of complex dimension $d+1$. There is a constant $c>0$ independent of $d$
such that for every $\gamma>0$ and every holomorphic map
$\Phi:X\to \Proj(V)$ of generic rank $n$ we can find a subset $E$ of
$\Proj(V^{*})$
satisfying the following properties:
\begin{enumerate}
\item[(a)] $\sigma_\FS(E)\leqslant c\,d^2\,e^{-\gamma/c}$\,;
\item[(b)] If $[\xi]$ is outside $E$, the current $\Phi^*[H_\xi]$ is well-defined and we have
$$\big\|\Phi^*[H_\xi]-\Phi^*(\omega_\FS)\big\|_{U,-2}\leqslant \gamma\,.$$
\end{enumerate}
\end{theorem}

The proof of the above result uses some properties of quasi-psh functions that we will recall here. For the details, see \cite{DS:06}. For simplicity, we will state these properties for $\Proj^d$ but we will use them for $\Proj(V^{*})$.

A function $u:\Proj^d\to \R\cup\{-\infty\}$ is called {\it quasi-psh} if it is locally the difference of a psh function and a smooth function. For such a function $u$, there is a constant $c>0$ such that $\ddc u+c\omega_{\FS}$ is a positive closed $(1,1)$-current.  Following Proposition A.3 and Corollary A.4 in \cite{DS:06} (in these results we can choose $\alpha=1$), we have:

\begin{proposition} \label{prop_qpsh}
There is a universal constant $c>0$ independent of $d$ such that for any quasi-psh function $u$ with $\max u=0$ and $\ddc u\geqslant -\omega_\FS$ we have
$$\|u\|_{L^1}\leqslant {\frac12}\big(1+\log d\big) \quad \mbox{and}\quad \|e^{-u}\|_{L^1}\leqslant cd\,,$$
where the norm $L^1$ is with respect to $\sigma_\FS$.
\end{proposition}

Quasi-psh functions are quasi-potentials of positive closed $(1,1)$-currents. For such a current, we will use the following notion of mass
\[
\|S\|:=\big(S,\omega_\FS^{d}\big)
\]
which is equivalent to the classical mass norm.
More precisely, we have the following result.

\begin{lemma} [$\partial\overline\partial$-Lemma for currents]\label{lemma_pot}
Let $S$ be a positive closed $(1,1)$-current on $\Proj^d$. Assume that the mass of $S$ is equal to $1$, that is, $S$ is cohomologous to $\omega_\FS$. Then there is a unique quasi-psh function $v$ such that
$$\max v=0 \quad\mbox{and}\quad \imat\partial\db v=S-\omega_\FS.$$
\end{lemma}

Finally, we will need the following lemma.

\begin{lemma} \label{lemma_key}
Let $\Sigma$ be a closed subset of $\Proj^d$ and $u\in L^1(\Proj^d)$ a function which is continuous on $\Proj^d\setminus \Sigma$. Let $\gamma$ be a positive constant. Suppose there is a positive closed $(1,1)$-current $S$ of mass $1$ such that $-S\leqslant \imat\partial\db u\leqslant S$ and
$\int_{\Proj^d} u\,d\sigma_\FS=0$.  Then, there is a universal constant $c>0$ and a Borel set $E\subset \Proj^d$ depending only on $S$ and $\gamma$ such that
$$\sigma_\FS(E)\leqslant cd^2e^{-\gamma}\quad \mbox{and}\quad |u(a)|\leqslant \gamma$$
for $a\not\in \Sigma\cup E$.
\end{lemma}
\proof
Let $v$ be as in Lemma \ref{lemma_pot}. Define $m:=\int_{\Proj^d} v\,d\sigma_\FS$. By Proposition \ref{prop_qpsh}, we have
$-{\frac12}(1+\log d)\leqslant m\leqslant 0$. Define $w:=u+v$. We have $\imat\partial\db w\geqslant -\omega_\FS$. Since $u$ is continuous outside $\Sigma$, the last property implies that $w$ is equal to a quasi-psh function outside $\Sigma$. We still denote this quasi-psh function by $w$. Define $l:=\max w$. Applying
Proposition \ref{prop_qpsh} to $w-l$, we obtain that
$$m-l=\int (w-l) d\sigma_\FS\geqslant - {\frac12}\big(1+\log d\big).$$
It follows that
$$l\leqslant {\frac12}\big(1+\log d\big).$$
We have
$$u=w-v\leqslant l-v\leqslant {\frac12}\big(1+\log d\big) -v.$$
Let $E$ denote the set $\{v<-\gamma+{\frac12}(1+\log d)\}$. This set depends only on $\gamma$ and on $S$. We have $u\leqslant \gamma$ outside $\Sigma\cup E$. The same property applied to $-u$ implies that
$|u|\leqslant \gamma$ outside $\Sigma\cup E$. It remains to bound the size of $E$.
The last estimate in Proposition \ref{prop_qpsh} yields
$$\sigma_\FS(E)\lesssim d \exp\big(-\gamma+{\tfrac12}(1+\log d)\big)\lesssim d^2e^{-\gamma}.$$
This is the desired inequality.
\endproof

\begin{proof}[Proof of Theorem \ref{th_equi}] Fix an open set $W$ such that $U\Subset W\Subset X$. Observe that when $\Phi^{-1}(H_\xi)$ does not contain any open subset of $X$ then $\Phi^*[H_\xi]$ is well-defined. Indeed, we can write locally $[H_\xi]=\imat\partial\db u$ for some psh function $u$ and define $\Phi^*[H_\xi]:=\imat\partial\db (u\circ\Phi)$. The function $u$ is smooth outside $H_\xi$ and equal to $-\infty$ on $H_\xi$. So, the expression $\imat\partial\db (u\circ\Phi)$ is meaningful since $u\circ\Phi$ is not identically $-\infty$, see \cite{Meo96} for details. Let $\Sigma$ denote the closure of the set of $[\xi]$ which do not satisfy the above condition. By shrinking $X$ we can assume that $\Sigma$ is an analytic set with boundary in $\Proj(V^{*})$. In particular, its volume is equal to 0. The currents $\Phi^*[H_\xi]$ depend continuously on $\xi\in \Proj(V^{*})\setminus \Sigma$.

Fix a smooth positive $(n,n)$-form $\nu$ with compact support in $W$ such that for any $\cC^2$ real form $\varphi$ of bidegree $(n-1,n-1)$ with compact support in $U$ and  $\|\varphi\|_{\cC^2}\leqslant 1$ we have $-\nu\leqslant \imat\partial\db\varphi\leqslant \nu$. Let $M$ be the analytic subset of points $(x,\xi)$ in $X\times \Proj(V^{*})$ such that $x\in H_\xi$. It is of dimension $n+d-1$.
Let $\pi_1$ and $\pi_2$ denote the natural projections from $M$ onto $X$ and $\Proj(V^{*})$ respectively. Define
$v:=(\pi_2)_*(\pi_1)^*(\varphi)$. This is a function on $\Proj(V^{*})$ whose value at $\xi\in \Proj(V^{*})\setminus\Sigma$ is the integration of $\pi_1^*(\varphi)$ on the fiber $\pi_2^{-1}(\xi)$. So, we have
$$v(\xi)=\big(\Phi^*[H_\xi],\varphi\big)\,.$$
Hence, $v$ is continuous on $\Proj(V^{*})\setminus\Sigma$. Since the form $\omega_\FS$ on $\Proj(V)$ is the average of $[H_\xi]$ with respect to the measure $\sigma_\FS$ on $\xi$, we also have
$$\int vd\sigma_\FS=\big( \Phi^*(\omega_\FS),\varphi\big)\,.$$
Define $T:=(\pi_2)_*(\pi_1)^*(\nu)$. Recall that $\nu$ is positive. It is closed since it is of maximal bidegree. It follows that $T$ is a positive closed $(1,1)$-current on $\Proj(V^{*})$, since the operators $(\pi_1)^*$ and $(\pi_2)_*$ preserve the positivity. Since $-\nu\leqslant\imat\partial\db\varphi\leqslant \nu$, we have $-T\leqslant \imat\partial\db v\leqslant T$. Let $m$ denote the mass of $\nu$ considered as a positive measure.
If $a$ is a generic point in $X$, then $T$ is cohomologous to $m(\pi_2)_*(\pi_1)^*(\delta_a)$ where $\delta_a$ is the Dirac mass $a$, because $\nu$ is cohomologous to $m\delta_{a}$. The last expression is $m$ times the current of integration of the hyperplane $H$ of points $\xi$ such that $\Phi(a)\in H_\xi$. So, the mass of $T$ is equal to $m$.  In particular, it is independent of $\Phi$ and $\varphi$.

Define the function $u$ on $\Proj(V^{*})$ by
\[
u:={\frac1m} \big(v-(\Phi^*(\omega_\FS),\varphi)\big)\,.
\]
We deduce from the above discussion that $u$ satisfies the hypothesis of Lemma \ref{lemma_key} for $S:={\frac1m}T$. Applying this lemma to $\gamma/m$ instead of $\gamma$, we find a set $E'$ independent of $\varphi$ such that
$\sigma_\FS(E')\leqslant cd^2e^{-\gamma/m}$ and $|u|\leqslant \gamma/m$ outside $\Sigma\cup E'$. It follows that
$$\big\|\Phi^*[H_\xi]-\Phi^*(\omega_\FS)\big\|_{U,-2}\leqslant \gamma$$
for $\xi$ out of $E:=\Sigma\cup E'$. It is enough to replace $c$ by $\max(c,m)$ in order to obtain the theorem.
\end{proof}
We show now how to apply Theorem \ref{th_equi} to prove Theorem \ref{conv-speed1}.
\begin{proof}[Proof of Theorem \ref{conv-speed1}]
We shall apply Theorem \ref{th_equi} for $V=H^0(X,L^N)^*$ and for the Kodaira map $\Phi_{N}:X\dashrightarrow\Proj H^0(X,L^N)^*$ defined in \eqref{kodmap} restricted to a suitable neighbourhood of $U$. Let $d_N=\dim H^0(X,L^N)^*$, hence, in the notation of Theorem \ref{th_equi}, $d+1=d_N$.
By Corollary \ref{kod-emb} (iii) we have 
\[
d_N\simeq N^n\,,\quad\log d_N\simeq \log N\,.
\]
For $[\xi]\in\mathbb{P}H^0_{(2)}(X,L^{N})$, the hyperplane $H_{\xi}\subset\mathbb PH^0_{(2)}(X,L^N)^*$ determines a current of integration $\Phi_{N}^{*}[H_{\xi}]$ on the zero-divisor of a section $s_{\xi}\in[\xi]$. This section is unique up to a multiplicative constant.

We know by Corollary \ref{kod-emb} (ii) that there exists a neighbourhood $W$ of $U$ and $N'(U)$ such that $\Phi_N:\;X\setminus Bs_{N}\longrightarrow \mathbb PH^0_{(2)}(X,L^N)^*$ is an embedding in $W$ for $N\geqslant N'(U)$. By \eqref{Tianest}, $N^{-1}\Phi_{N}^*(\omega_\FS)$ differs from $\frac{\imat}{2\pi}R^{L}$ on $U$ by a form of norm bounded by $C_{2,U}/N$. We apply
Theorem \ref{th_equi} for $\Phi_{N}|_W$ with $N\geqslant N'(U)$, $d=d_N-1$ and $\gamma=\lambda_N/2$. Thus, for $N\geqslant N'(U)$ there exist $E_N\subset\Proj H^0_{(2)}(X,L^N)$ such that
$\sigma_\FS(E_N)\leqslant c\,N^{2n}\,e^{-\lambda_{N}/c}$ and for all $[\xi_N]\in\Proj H^0_{(2)}(X,L^N)\setminus E_{N}$ we have
\[
\Big|\Big(\frac1N\Phi_{N}^*[H_{\xi_{N}}]-\frac1N\Phi_{N}^*(\omega_\FS),\varphi\Big)\Big|\leqslant \frac{\lambda_{N}}{2N}\|\varphi\|_{\cC^{2}}\,,\quad \varphi\in\Omega_0^{n-1,n-1}(U)\,.
\]
Hence
\[
\Big|\Big(\frac1N\Phi_{N}^*[H_{\xi_{N}}]-\frac{\imat}{2\pi}R^{L},\varphi\Big)\Big|\leqslant
\Big(\frac{C_{2,U}}{N}+\frac{\lambda_{N}}{2N}\Big)\|\varphi\|_{\cC^{2}}\,,\quad \varphi\in\Omega_0^{n-1,n-1}(U)\,.
\]
Choose now $N''(U)$ such that $\lambda_{N}\geqslant 2C_{2,U}$ for all $N\geqslant N''(U)$.
We obtain the items (a) and (b) of Theorem \ref{conv-speed1} by setting $N(U)=\max\{N'(U),N''(U)\}$.
 
Since
$\sum \sigma_\FS(E_N)<\infty$, the last property holds for a generic sequence $([\xi_N])\in\Omega$.
\end{proof}
\noindent
Note that by Remark \ref{rem-tianest} we can replace $C_{2,U}/N$ by $C_{2,U}/N^{2}$
if $\Theta=\frac{\imat}{2\pi}R^{L}$.
\begin{remark}
Let $1\leqslant p\leqslant n$ be an integer. Suppose there is a positive, closed $(n-p+1,n-p+1)$-form $\nu$ with compact support in $X$, which is strictly positive on $\overline{U}$. Then we can extend Theorem \ref{conv-speed1} to projective subspaces of codimension $p$ instead of hyperplanes $H_{\xi}$. This can be applied to obtain the equidistribution on $U$ of common zeros of $p$ random holomorphic sections i.e.\, of currents of the form $[s^{(1)}_{N}=\ldots=s^{(p)}_{N}=0]$.
\end{remark}
\begin{proof}[Proof of Theorem \ref{conv}]
Take an exhaustion $(U_{j})_{j\in\N}$ of $X$ by open relatively compact sets. By Theorem \ref{conv-speed1} there exist sets $\mathcal N_{j}\subset\Omega$, $j\in\N$, of $\mu$-measure zero such that \eqref{conveq} holds for all $\varphi\in\Omega_{0}^{n-1,n-1}(U_{j})$ and all sequences $\mathbf s\in \Omega\setminus \mathcal N_{j}$.
Now, since $\Omega_{0}^{n-1,n-1}(X)=\cup_{j\in\N} \Omega_{0}^{n-1,n-1}(U_{j})$, \eqref{conveq} holds for all $\mathbf s \in \Omega\setminus \mathcal N$ where $\mathcal N = \cup_{j\in\N} \mathcal N_{j}$.
\end{proof}
\begin{remark}\label{conv-var}
By using Remark \ref{rem-forms} (ii), (iii) we can obtain various versions of Theorems \ref{conv}, \ref{conv-speed1}.
Let $(X,\Theta)$ be an $n$-dimensional complete Hermitian manifold. Let $(L,h^L)$ be a Hermitian holomorphic line bundle over $X$.
\\[2pt]
(i) Assume that $(X,\Theta)$ and $(L,h^L)$ satisfy \eqref{i_bis}. 
Furthermore, assume that the spaces $H^0_{(2)}(X,L^N\otimes K_X)$ have finite dimension for all $N$ and $\dim H^0_{(2)}(X,L^N\otimes K_X)=O(N^n)$, $N\to\infty$. Then
the zero-divisors of generic random sequences $(s_{N}) \in \prod_{N=1}^\infty\Proj H^0_{(2)}(X,L^N\otimes K_X)$
are equidistributed with respect to $\tfrac{\imat}{2\pi}R^L$.
\\[2pt]
(ii) Assume that $L=K_X$, $h^{K_X}$ is induced by $\Theta$ and \eqref{i_bis} holds. 
Furthermore, assume that the spaces $H^0_{(2)}(X,K_X^N)$ have finite dimension for all $N$ and $\dim H^0_{(2)}(X,K_X^N)=O(N^n)$, $N\to\infty$. Then
the zero-divisors of generic random sequences $(s_{N}) \in \prod_{N=1}^\infty\Proj H^0_{(2)}(X,K_X^N)$
are equidistributed with respect to $\tfrac{\imat}{2\pi}R^{K_X}$.

In addition, the statement of Theorem \ref{conv-speed1} holds also in the situations (i) and (ii) above.

\end{remark}
\section{Zeros of pluricanonical sections on pseudoconcave and arithmetic quotients}\label{sec:psconc}
We recall the definition of pseudoconcavity in the sense of Andreotti and Grauert.
\begin{definition}[{Andreotti-Grauert \cite{AG:62}}]\label{def-ag}
Let $X$ be a complex manifold of complex dimension $n$ and $1\leqslant q\leqslant n$.
$X$ is called \emph{$q$-concave}
if there exists a smooth function
$\varphi:~X\longrightarrow (a,b]$, where $a\in\R\cup\{-\infty\}$, $b\in\R$,
such that for all $c\in]a,b]$ the superlevel sets
$X_c=\{\varphi>c\}$ are relatively compact in $X$, and
$\partial\db\varphi$ has at least $n-q+1$ positive eigenvalues
outside a compact set.
\end{definition}
\begin{example}\label{exconc}
Let $Y$ be a compact complex space and let $A\subset Y$ be an analytic subset of dimension $q$
which contains the singular locus of $Y$. Then $Y\setminus A$ is a $(q+1)$-concave manifold (see e.g.\ \cite[Prop.\,9]{Vaj94}).
\end{example}

We can formulate now the following consequence of Theorem \ref{conv}. We first recall some terminology.
Let $(X,J,\omega)$ be a Hermitian manifold, let $g^{TX}$ be the Riemannian metric associated to $\omega$ by \eqref{riem-met} and let $\ric$ the Ricci curvature of $g^{TX}$. The Ricci form $\ric_\omega$ is defined as the $(1,1)$-form associated to $\ric$ by
\begin{equation}\label{ric0}
\ric_\omega(u,v)=\ric(Ju,v)\,,\quad\text{for any $u,v\in T_{x}X$, $x\in X$}.
\end{equation}
If the metric $g^{TX}$ is K\"ahler, then
\begin{equation}\label{ric1}
\ric_\omega=\imat R^{K_{X}^{*}}=-\imat R^{K_{X}}
\end{equation}
where $K_{X}^{*}$ and $K_{X}$ are endowed with the metric induced by $g^{TX}$ (see e.g.\ \cite[Prob.\,1.7]{MM07}).
The metric $g^{TX}$, or the associated K\"ahler form $\omega$, are called \emph{K\"ahler-Einstein} if there exists a real constant $k$ such that
\begin{align*}
 &\ric=k\,g^{TX}  \\ \intertext{or, equivalently}
 &\ric_\omega=k\,\omega\,.
\end{align*}
\begin{theorem}\label{ke}
Let $(X,\omega)$ be an $n$-dimensional complete K\"ahler manifold with $\ric_\omega\leqslant -k\,\omega$ for some constant $k>0$.  Assume that $X$ is $(n-1)$-concave. Then the zero-divisors of generic random sequences $(s_{N}) \in \prod_{N=1}^\infty \Proj H^0_{(2)}(X,K_{X}^N)$
are equidistributed with respect to $-\frac{1}{2\pi}\ric_\omega$.
If $(X,\omega)$ is K\"ahler-Einstein with $\ric_\omega=-k\,\omega$, $k>0$, then we have equidistribution with respect to $\frac{k}{2\pi}\omega$\/.
\end{theorem}
The $L^{2}$ inner product here is constructed with respect to the volume form of $\omega$
and the Hermitian metric on $K_{X}$ induced by $\omega$.
\begin{proof}
The proof follows immediately from Remark \ref{conv-var} (ii), since $\imat R^{K_{X}}=-\ric_\omega$\,.
Moreover, an important feature of any $q$-concave manifold $X$ is that for every holomorphic line bundle $F$ on $X$ the space of holomorphic sections $H^{0}(X,F)$ is of finite dimension. This is a consequence of the finiteness theorem of Andreotti-Grauert \cite{AG:62}.  In \cite[Th.\,3.4.5]{MM07} this is shown by elementary means (going back to Andreotti) which also deliver the cohomology growth condition \eqref{upb}, cf.\ \cite[(3.4.3)]{MM07}.
\end{proof}
\begin{example}\label{kdample}
Assume that $X$ is an $n$-dimensional Zariski open set in a compact complex space $X^{*}$ such that the codimension of the analytic set $X^{*}\setminus X$ is at least two. Then $X$ is $(n-1)$-concave. Assume moreover that (e.g.\ after desingularization of $X^{*}$) $X$ is biholomorphic to $\widetilde{X}\setminus D$, where $\widetilde{X}$ is a compact manifold, $D$ is an effective divisor with only normal crossings such that $K_{\widetilde{X}}\otimes\cO_{\widetilde{X}}(D)$ is ample. Then $X$ admits a (unique up to constant multiple) K\"ahler-Einstein metric with negative Ricci curvature, by a theorem due to R. Kobayashi \cite{KobR84} and Tian-Yau \cite{TiaYau87}.
\end{example}
\begin{example}
Another class of examples is the following. By a well-known theorem of Cheng-Yau \cite{CY80} and Mok-Yau \cite{MY83}, any Riemann domain
$\pi:D\to\C^{n}$ with $\pi(D)$ bounded (e.g.\ a bounded domain of holomorphy in $\C^{n}$) carries a unique complete K\"ahler-Einstein metric $\omega$ of constant negative Ricci curvature, $\ric_\omega=-\omega$\,.

Assume that $X$ is an $n$-dimensional Zariski open set in a compact complex space $X^{*}$ such that the codimension of the analytic set $X^{*}\setminus X$ is at least two and $X$ is covered by a bounded domain of holomorphy in $\C^{n}$.
Then $X$ is $(n-1)$-concave and it carries a complete K\"ahler-Einstein metric $\omega$ of constant negative Ricci curvature, $\ric_\omega=-\omega$\,.

\end{example}
Other examples where Theorem \ref{ke} applies are provided by arithmetic quotients.
Let us quote the following fundamental result about compactification of arithmetic quotients.
\begin{theorem}[{Satake \cite{Sa:60}, Baily-Borel \cite{BB:66}}]
Let $X$ be an $n$-dimensional, $n\geqslant2$, irreducible  quotient of a bounded symmetric domain
$D$ by a torsion-free arithmetic group $\Gamma\subset\aut(D)$. Then there exists a compactification $X^*$ of $X$ such that: (i) $X^*$ is a normal projective variety and (ii) $X^*\setminus X$ is a complex analytic variety  is of codimension $\geqslant2$ in $X^*$.
\end{theorem}
We are now in the position to prove Corollary \ref{equi-arit}.
\begin{proof}[Proof of Corollary \ref{equi-arit}]
By Example \ref{exconc} it follows that an $n$-dimensional arithmetic quotient $X$, with $n\geqslant2$,
is $(n-1)$-concave, since the singular locus of $X^*$ has dimension at most $n-2$. (The proof of pseudoconcavity (initially in a weaker sense) of arithmetic quotients by Andreotti-Grauert \cite{AG:61} and  Borel \cite{Bo:70} actually also shows that $X$ is $(n-1)$-concave.)

Assume that $X=\Gamma\backslash  D$, where $D$ is a bounded symmetric domain and $\Gamma$ an arithmetically defined discontinuous group that acts freely on $D$. The domain
$D$ has a canonical $\aut(D)$-invariant K\"ahler metric, namely the Bergman metric.
Pick a basis $\{S_{i}\}_{i\geqslant1}$ of the Hilbert space $H^0_{(2)}(D)$ of square-integrable holomorphic functions. The Bergman kernel is defined as the locally uniformly convergent sum
\[ P(z,w)=\sum {S_{i}}(z)\,\overline{S}_{i}(w)\] which is in fact independent of the choice of the basis. The Bergman metric
\begin{equation}\label{bergmet}
\omega^{\mathcal B}_{D}:=\imat\partial\bar{\partial} \log P(z,z)
\end{equation}
is invariant under $\aut(D)$ (\cite[Ch.\,4, \S1, Prop.\,2]{Mok89b}) and obviously K\"ahler. Therefore it descends to a K\"ahler metric on $X$, denoted $\omega^{\mathcal B}_{X}$. Since $(D,\omega^{\mathcal B}_{D})$ is a Hermitian symmetric manifold, the Bergman metric is complete. Therefore $(X,\omega^{\mathcal B}_{X})$ is complete with respect to the Bergman metric.

The volume form of $\omega^{\mathcal B}_{D}$ defines a Hermitian metric in the canonical line bundle $K_{D}$ and there holds $\imat R^{K_{D}}= \omega^{\mathcal B}_{D}$.
In fact, since $\omega^{\mathcal B}_{D}$ is invariant, the forms $(\omega^{\mathcal B}_{D})^n$ and $P(z,z)dz_{1}\wedge d\bar{z}_{1}\wedge \dots \wedge dz_{n}\wedge d\bar{z}_{n}$ are both invariant $(n,n)$ forms (\cite[Ch.\,4, \S1, Prop.\,2]{Mok89b}). Since $D$ is homogeneous, an invariant object is determined by its value in a point. Hence, two $(n,n)-$forms only differ by a constant. Writing
\[
(\omega^{\mathcal B}_{D})^n =c\,P(z,z)dz_{1}\wedge d\bar{z}_{1}\wedge \dots \wedge dz_{n}\wedge d\bar{z}_{n}
\]
shows that $\imat R^{K_{D}}= \imat\partial\bar{\partial} \log P(z,z)=\omega^{\mathcal B}_{D}$ (\cite[Ch.\,4, \S 1, Prop.\,3]{Mok89b}). Since $\omega^{\mathcal B}_{D}$ descends to $X$, we also have $\imat R^{K_{X}}=\omega^{\mathcal B}_{X}$. The result follows therefore from Theorems
\ref{conv}, \ref{conv-speed1} and \ref{ke}.
\end{proof}
Note that in the above examples the base manifold $X$ turns out to be actually quasi-projective.
In fact, we can replace the hypothesis that $X$ is $(n-1)$-concave in Theorem \ref{ke} by the hypothesis that $X$ is quasi-projective and obtain the same conclusions. More precisely we have the following.
\begin{theorem}\label{ke-proj}
Let $(X,\omega)$ be an $n$-dimensional complete K\"ahler manifold with $\ric_\omega\leqslant -k\,\omega$ for some constant $k>0$.  Assume that $X$ is is quasi-projective. Then the zero-divisors of generic random sequences $(s_{N}) \in \prod_{N=1}^\infty \Proj H^0_{(2)}(X,K_{X}^N)$
are equidistributed with respect to $-\frac{1}{2\pi}\ric_\omega$.
If $(X,\omega)$ is K\"ahler-Einstein with $\ric_\omega=-k\,\omega$, $k>0$, then we have equidistribution with respect to $\frac{k}{2\pi}\omega$\/. Moreover, we have an estimate of the convergence speed on compact sets as in Theorem \ref{conv-speed1}.
\end{theorem}
\begin{proof}
Let $\overline{X}$ be a smooth compactification of $X$ such that $D=\overline{X}\setminus X$
is a divisor with simple normal crossings.
By \cite[Lemma\,5.1]{Tia:90} any holomorphic section of $H^0_{(2)}(X,K_{X}^N)$ extends to
a meromorphic section of $K_{\overline{X}}^N$ with poles along of order at most $N$ along $D$, i.e., $H^0_{(2)}(X,K_{X}^N)\subset H^0(\overline{X},K_{\overline{X}}^N\otimes[D]^N)$. It follows that $H^0_{(2)}(X,K_{X}^N)$ are finite dimensional and \eqref{upb} is satisfied. By applying Theorems
\ref{conv} and \ref{conv-speed1} (cf.\ Remark \ref{conv-var}) we obtain the result.
\end{proof}
\begin{remark}
We can identify more precisely the sections of $H^0_{(2)}(X,K_{X}^N)$ in Theorem \ref{ke-proj} in terms of sections on a compactification $\overline{X}$ as above.
\\[2pt]
(i) By \cite[Prop.\,1.11]{Nad:90} we have $H^0_{(2)}(X,K_{X}^N)\subset H^0(\overline{X},K_{\overline{X}}^N\otimes[D]^{N-1})$ in the conditions of Theorem \ref{ke-proj}.
\\[2pt]
(ii) Let $\overline{X}$ be a projective manifold and let $D$ be an effective divisor with only normal crossings such that $K_{\overline{X}}\otimes\cO_{\overline{X}}(D)$ is ample. Let $\omega$ be a complete K\"ahler metric on $X=\overline{X}\setminus D$ such that $\ric_\omega=-\omega$ (cf.\ \cite{KobR84,TiaYau87}). 
The metric $\omega$ has Poincar\'e growth on $X$ (cf.\ the proof of \cite[Th.\,1]{KobR84}) so by the proof of \cite[Lemma\,5.2]{Fu:92} we have $H^0_{(2)}(X,K_{X}^N)= H^0(\overline{X},K_{\overline{X}}^N\otimes[D]^{N-1})$.
\end{remark}
\section{Equidistribution of zeros of modular forms}\label{sec:modular}

Consider the group $SL_2(\mathbb Z)$ acting on the hyperbolic plane $\mathbb H$ via linear fractional transformations. Let $\Gamma$ be any subgroup of finite index that acts freely and properly discontinuously. Then the quotient space is naturally endowed with a smooth manifold structure and we want to consider holomorphic line bundles on $Y=\Gamma \backslash \mathbb H.$\par
Denote by $p$ the projection map $ \mathbb H \to \Gamma \backslash \mathbb H$.  For any line bundle $L\rightarrow Y$ there exists a global trivialization \[\varphi: p^*L \to \mathbb H \times \mathbb C.\] By invariance, $(p^*L)_{\tau}= (p^*L)_{\gamma\tau}$ for any $\gamma \in \Gamma$ and we can form $j_{\gamma}(\tau):=\varphi_{\gamma\tau}\circ\varphi_{\tau}^{-1} .$ Clearly, $j_{\gamma} \in \mathcal O^*(\mathbb H)$ and satisfies $ j_{\gamma\gamma^{\prime}}(\tau)= j_{\gamma}(\gamma^{\prime}\tau)j_{\gamma^{\prime}}(\tau).$
The map $j:\Gamma\times \mathbb H  \to \mathbb C^{\times}$ is called an automorphy factor for $\Gamma$. Conversely, any automorphy factor induces a $\Gamma-$ action on the trivial bundle $\mathbb H \times \mathbb C$ via $(\tau,z)\mapsto (\gamma\tau,j_{\gamma}(\tau)z)$ and the quotient becomes a holomorphic line bundle on $Y$ with transition functions given by the $j_{\gamma}.$\par
Holomorphic sections of $L$ can be identified with holomorphic functions on $\mathbb H$ that satisfy $f(\gamma \tau)= j_{\gamma}(\tau)f(\tau).$ Sections of the tensor powers $L^{N}$ satisfy $f(\gamma \tau)= j_{\gamma}(\tau)^{ N}f(\tau).$ Hermitian metrics on $L$ are identified with smooth (real) functions that satisfy $h(\gamma \tau)=| j_{\gamma}(\tau)|^{-2}h(\tau), $ the induced metric on $L^{N}$ corresponds to $h^N(\gamma \tau)=| j_{\gamma}(\tau)|^{-2N}h^N(\tau). $\par
We consider a canonical map $\Gamma\times\mathbb H \to \mathbb C^{\times}$,\[(\gamma,\tau)\mapsto \left.\frac{d\gamma}{dz}\right|_{\tau}\]  where $\frac{d\gamma}{dz}$ is the complex differential. The chain rule implies $ \frac{d\gamma\gamma^{\prime}}{dz}(\tau)=\frac{d\gamma}{dz}(\gamma^{\prime}\tau) \frac{d\gamma^{\prime}}{dz}(\tau).$ We call this map the canonical automorphy factor. Recall that $\gamma \in \Gamma$ are the transition maps of the coordinate charts on $Y$. Therefore the $\frac{d\gamma}{dz}$ are the transition functions of the tangent bundle $TY$. Explicitly,
\[ \gamma \tau=
\left( \begin{array}{cc}
    a &b \\
    c& d \\
  \end{array}\right)\tau=
\frac{a\tau+b}{c\tau+d}\,,
\qquad \frac{d\gamma}{dz}(\tau)=\frac{1}{(c\tau+d)^{2}}\,\cdot\]
The transition functions of the dual bundle $T^*Y=K_{Y}$ are then $j_{\gamma}(\tau)= (c\tau+d)^{2}$.

Summarizing, holomorphic functions on $\mathbb H$ that satisfy the transformation law $f(\gamma \tau)= (c\tau+d)^{2N}f(\tau)$ are in correspondence with holomorphic sections of the bundle $K^{\otimes N}_{Y}\to Y.$
\begin{definition}\label{def-modular}
 A modular form for $\Gamma$ of weight $2N$ is a function on $\mathbb H$ such that:
\\[2pt]
\ (1) $f(\gamma \tau)= (c\tau+d)^{2N}f(\tau)$ , all $\tau \in \mathbb H$,
 \\[2pt]
(2) $f$ is holomorphic on $\mathbb H$,
 \\[2pt]
(3) $f$ is ``holomorphic at the cusps''.
\\[2pt]
A modular form $f$ that is zero at every cusp of $\Gamma$ is called a cusp form. We write $\mathcal M_{2N}(\Gamma)$ for the space of modular forms of weight $2N$ and $\mathcal S_{2N}(\Gamma)$ for the subspace of cusp forms.
\end{definition}
The last condition means the following. If the index $[SL_{2}(\mathbb Z):\Gamma]$ is finite, then for some natural number $\ell$ the transformation $z \mapsto z+\ell$ is contained in
$\Gamma$ and by (1), $f(z+\ell)=f(z).$ We identify a neighbourhood $\{ 0<\Re z \leqslant \ell, \Im z > c\}$ of $\infty$  with the punctured disk $\{ 0<|w|<\exp(-2\pi c/\ell)\}$ via the map $q: z \mapsto \exp(2\pi \imat z/\ell)$ and define the $q$-expansion $\hat{f}$ of $f$ at the cusp $\infty$ by $q^*\hat{f}=f$. By definition, $f$ is holomorphic (resp. $0$) at $\infty$ if $\hat{f}$ is holomorphic (resp. $0$) at $0.$ Now any other cusp $\sigma$ of $\Gamma$ is of the form $\sigma = \alpha\infty$ with some $\alpha \in SL_{2}(\mathbb Z)$ and we say that $f$ is holomorphic (resp.\ $0$) at $\sigma$ if $j_{\alpha}^{-1}f\circ\alpha$ is is holomorphic (resp.\ $0$) at $\infty$.

The space $\mathcal M_{2N}(\Gamma)$ is finite dimensional for any $N$, in case $N=0$ the dimension is 1 and in case $N<0$ the dimension is $0$. In the following, we consider $N\geqslant 1$.\par
The bundle $K_{Y}$ inherits a positively curved metric from the hyperbolic space. Namely, we endow the canonical bundle of $\mathbb H$ (which is trivial) with the metric $h$ described in terms of the length of the section $1$ by $|1|_{h}(z):=|\Im z|=|y|$. This metric descends to a metric on $K_{Y}$ and we have
\begin{align*}
p^{*}(\imat R^{K_{Y}}) &= -\imat\partial\bar{\partial} \log y^2\\
 &= -\imat\partial\bar{\partial} \log (-\frac{1}{4}(z-\bar z)^2)\\
 &=-\imat\frac{2}{(z-\bar{z})^2}dz\wedge d\bar z\\
 &=\imat\frac{1}{2y^2}dz\wedge d\bar z\\
 &=\frac{dx\wedge dy}{y^2}\;\cdot
 \end{align*}
Let $D$ be a fundamental domain for $\Gamma.$ If $f$ is any modular form of weight $2N$, $g$ a cusp form of the same degree, then the integral 
\[
\int_{D} f(z)\bar{g}(z)y^{2N-2}dxdy
\] 
converges and defines a Hermitian product on $\mathcal S_{2N}(\Gamma)$, called the Petersson inner product. This is just the induced $L^2$ product on $H^0(Y,K^{N}_{Y}).$  We wish to describe the subspace $H_{(2)}^0(Y,K^{N}_{Y})$ of square integrable sections, i.e. $f\in H_{(2)}^0(Y,K^{N}_{Y})$ is a function on $\mathbb H$ satisfying conditions (1),(2) from Definition \ref{def-modular} and such that
\[\int_{D}|f|^2y^{2N-2}dxdy < \infty.\] Since the cusps of $\Gamma$ are all of the form $\sigma_{i}\infty$ with $\sigma_{i}\in SL_{2}(\mathbb Z), 1\leqslant i\leqslant |SL_{2}(\mathbb Z):\Gamma|,$ it suffices to consider the integral in a neighbourhood $U(\infty)$ of $\infty.$ Using $\hat{f}$ defined on $A:=\{ 0<|w|<\exp(-2\pi c/\ell)\}$ we compute
\begin{align*}
  \int_{U(\infty)\cap D}|f|^2y^{2N-2}dxdy  &\sim \int_{c}^{\infty}\int_{0}^{\ell}|f|^2y^{2N-2}dxdy \\
  &\geqslant \int_{q^{-1}(A)}|f|^2dx\wedge dy \qquad \qquad \text{(since  $N\geqslant1$)}\\
  & =\int_{A}|\hat{f}(w)|^2\frac{\imat}{2}\frac{dw\wedge d\bar{w}}{4\pi^2|w|^2}\\
  &\gtrsim  \int_{A}{|\widehat{f}(w)|^2}{\frac{\imat}{2}dw\wedge d\bar{w}}\,.
\end{align*}
Hence, at each cusp  the $q$-expansion of $f$ is locally square integrable. We have the following simple lemmas.
\begin{lemma}
 Denote $\D = \{z\in \mathbb C: |z|<1\}$ and let $\widehat{f}$ be holomorphic in $\D\setminus \{0\}.$ Then $\widehat{f}\in L^2_{loc}(\D)$ if and only if $\widehat{f}$ can be holomorphically extended to the whole disk.
\end{lemma}
\begin{proof}
 Consider the Laurent-expansion $\widehat{f}= \sum_{\nu \in \mathbb Z}a_{\nu}w^{\nu}$. In the annuli $R_{r_{1},r_{2}}= \{r_{1} \leqslant |w| \leqslant r_{2}\}$ the series is normally convergent. We compute
 \begin{align*}
 \int_{R_{r_{1},r_{2}}} |\widehat{f}|^2 \tfrac{\imat}{2}dw\wedge d\bar{w}&= \sum_{\mu,\nu}a_{\mu}\bar{a}_{\nu}\Big(\int w^{\mu}\bar{w}^{\nu} \tfrac{\imat}{2}dw\wedge d\bar{w}\Big)\\
 &= \sum_{\nu \in \mathbb Z} |a_{\nu}|^2 2\pi \int_{r_{1}}^{r_{2}} r^{2\nu+1}dr.
 \end{align*}
 Either the principal part of $\widehat{f}$ at $0$ is zero or there exists $\nu <0$ such that $|a_{\nu}|^2> 0.$ In this case, \[\| \widehat{f}\|^2_{L^2(\D)}\geqslant 2\pi |a_{\nu}|^2\int_{r_{1}}^{r_{2}} r^{2\nu+1}dr\]and this is unbounded as $r_{1}\to 0$.
\end{proof}
Thus $H^0_{(2)}(Y,K^{N}_{Y})$ is of finite dimension. More precisely there holds:
\begin{lemma}
\label{cusps}
With the above notations we have $H^0_{(2)}(Y,K^{N}_{Y}) \simeq \mathcal{S}_{2N}(\Gamma)$.
\end{lemma}
\begin{proof}
Since the Petersson inner product is the $L^2$ product we have $S_{2N}(\Gamma)\subset H^0_{(2)}(Y,K^{N}_{Y})$. Let $f \in H^0_{(2)}(Y,K^{N}_{Y})$. By the preceding lemma, $f$ corresponds to a modular form of weight $2N$. We have to show that this modular forms vanishes at the cusps. It suffices to consider the cusp at $\infty$. Assume that $f$ does not vanish at $\infty$ and $\widehat{f}$ is the $q-$expansion of $f$ around $\infty$. Then $|\widehat{f}(0)|^2=b>0$ and $|f|^2\geqslant b/2$ in a neighbourhood of $\infty$. Then \[ \int_{c}^{\infty}\int_{0}^\ell |f|^2 y^{2N-2}dxdy \geqslant \frac{\ell b}{2}\int_{c}^{\infty}y^{2N-2}dy=\infty,\quad N\geqslant 1.\]
This contradiction shows that $f$ corresponds to a cusp form.
\end{proof}
It is well-known, by the Riemann-Roch formula, that $\dim\mathcal{S}_{2N}(\Gamma)\sim N$, as $N\to\infty$.
Therefore, Theorems \ref{conv}, \ref{conv-speed1} (cf.\ Remark \ref{conv-var}) and the above discussion immediately imply Corollary \ref{equi-modular}.
%
\section{Equidistribution on quasiprojective manifolds}\label{sec:quasiproj}
In the previous sections we considered the equidistribution with respect to some canonical
K\"ahler metrics. We turn now to the case of quasiprojective manifolds and construct adapted metrics using a method of Cornalba and Griffiths.
They depend of some choices but have the advantage of being very general.

Let $X\subset\Proj^{k}$ be a quasiprojective manifold, denote by $\overline{X}\subset\Proj^{k}$
its projective closure and let $\Sigma=\overline{X}\setminus X$. Denote by $L=\mathcal{O}(1)|_{\overline{X}}$ the restriction of the hyperplane line bundle $\mathcal{O}(1)\to\Proj^{k}$.

We consider a resolution
of singularities $\pi:\widetilde{X}\to\overline{X}$ in order to construct appropriate metrics on $X$ and $L|_{X}$. More precisely, there exists a finite
sequence of blow-ups
\[\widetilde{X}= X_{m} \stackrel{\tau_{m}}{\longrightarrow} X_{m-1}\stackrel{\tau_{m-1}}{\longrightarrow} \cdots \stackrel{\tau_{1}}{\longrightarrow}X_{0}=\overline{X}\]
along smooth centers $Y_{j}$ such that
\begin{enumerate}
\item $Y_{j}$ is contained in the strict transform $\Sigma_{j}=\overline{\tau_{j}^{-1} (\Sigma_{j-1}\backslash Y_{j-1})}$, $\Sigma_0=\Sigma$,
\item the strict transform of $\Sigma$ through $\pi = \tau_{m}\circ\tau_{m-1}\circ \cdots \tau_{1}$ is smooth and $\pi^{-1}(\Sigma)$ is a divisor with simple
normal crossings in $\widetilde{X}$,
\item $X=\overline{X}\backslash \Sigma \simeq \widetilde{X}\backslash\pi^{-1}(\Sigma)$ are biholomorphic.
\end{enumerate}

If $\pi^{-1}(\Sigma)= \cup S_{j}$ is a decomposition into smooth irreducible components, we denote for each $j$ the associated holomorphic line
bundle by $\mathcal O_{\widetilde{X}}(S_{j})$ and by $\sigma_{j}$ a holomorphic section vanishing to first order along $S_{j}$. Let $\Theta$ be the
fundamental form of any smooth Hermitian metric on $\widetilde X$.
Since $\widetilde X$ is projective, we can choose a K\"ahler form $\Theta$, but the construction works in general.
The generalized Poincar\'{e} metric on $\widetilde{X}\backslash\pi^{-1}(\Sigma)\simeq X$ is defined by the Hermitian form
\begin{equation}
\label{gpm}
\Theta_{\varepsilon} = \Theta - \imat\,\varepsilon\sum_{j}\partial\overline{\partial}\log(-\log\| \sigma_{j}\|^2_{j} )^2\,,\quad \text{$0<\varepsilon\ll1$ fixed},
\end{equation}
where we have chosen smooth Hermitian metrics $\|\cdot\|_{j}$ on $\mathcal O_{\widetilde{X}}(S_{j})$ such that $\| \sigma_{j}\|_{j}<1$. The generalized Poincar\'{e} metric \eqref{gpm} is a  complete Hermitian metric on $\widetilde{X}\backslash\pi^{-1}(\Sigma)\simeq X$  and satisfies the curvature estimates
 \[ -C\Theta_{\varepsilon} < \imat R^{K_{X}} <C\Theta_{\varepsilon} , \qquad |\partial \Theta_{\varepsilon}|_{\Theta_{\varepsilon}}<C \]
with some positive constant $C$ (where the metric on $K_{X}$ is the induced metric by $\Theta_{\varepsilon}$). A proof of this fact can be found in \cite[Lemma\,6.2.1]{MM07}.
Next we construct a metric on $L$ that dominates the Poincar\'{e} metric. By \cite[Lemma\,6.2.2]{MM07}, there exists a Hermitian line bundle $(\widetilde L, h^{\widetilde L})$ on $\widetilde X$ with positive curvature $R^{\widetilde L}$ on $\widetilde X$, and such that
\[
\widetilde L|_{\widetilde{X}\backslash \pi^{-1}(\Sigma)} \simeq \pi^*(L^m)|_{X}\,,
\]
with some $m\in \mathbb N$. If we equip $L|_{X}$ with the metric
\begin{equation}\label{qp1.1}
h_{\delta}^L = (h^{\widetilde L})^{\frac{1}{m}}\prod_{j} (-\log\|\sigma_{j}\|_{j}^2)^{2\delta}, \qquad \text{for $0<\delta\ll1$},
\end{equation}
we obtain that for any $0<\delta\ll1$ there exists $\eta>0$ such that following estimate holds
\[ \imat R^{h_{\delta}^L}>\eta\Theta_{\varepsilon}\,.\]
This holds true since $R^{\widetilde L}$ extends to a strictly positive $(1,1)$-form dominating a small positive multiple of $\Theta$ on $\widetilde X$.
Thus, the expansion of the Bergman kernel holds as in Theorem \ref{bke}. Moreover, the space of holomorphic $L^2$-sections
\[H^{0}_{(2)}(X,L^N,\Theta_{\varepsilon},h_{\delta}^{L^{N}}):= \big\{ s\in \mathcal O_{X}(L^N) \colon \int_{X} |s|_{h^{L^N}_{\delta}}^2 \Theta_{\varepsilon}^n/n!<\infty \big\} \]
is finite dimensional since the holomorphic $L^2$-sections extend holomorphically to all of $\widetilde{X}$, more precisely, $H^{0}_{(2)}(X,L^N)\subset H^{0}(\widetilde{X},\pi^{*}{L}^N)$, see \cite[(6.2.7.)]{MM07}.
In view of the previous discussion, Theorem \ref{conv} yields the following.
\begin{corollary}
Let $X\subset\Proj^{k}$ be a quasiprojective manifold. Denote by $L=\mathcal{O}(1)|_{X}$ the restriction of the hyperplane line bundle $\mathcal{O}(1)\to\Proj^{k}$. Fix metrics $\Theta_{\varepsilon}$ as in \eqref{gpm} and $h^{L}_{\delta}$ as in \eqref{qp1.1}. Then
the zero-divisors of generic random sequences $(s_{N}) \in \prod_{N=1}^\infty\Proj H^0_{(2)}(X,L^N,\Theta_{\varepsilon}, h_{\delta}^{L^{N}})$
are equidistributed with respect to $\frac{\imat}{2\pi}R^{(L,h_{\delta}^L)}$, where $\Theta_{\varepsilon}$ is the generalized Poincar\'e metric \eqref{gpm} and $h_{\delta}^L$ is defined by \eqref{qp1.1}.
\end{corollary}
\section{Equidistribution of zeros of orthogonal polynomials}\label{sec:poly}
We wish to illustrate the result of the previous Section in the case of polynomials. As mentioned in the Introduction the distribution of zeros of random polynomials is a classical subject.
Recent results were obtained by Bloom-Shiffman \cite{BL07} (see also \cite{Ber10}) concerning the equilibrium measure $\mu_{\rm eq}$ of a compact set $K$ endowed with a measure $\mu$ satisfying the Bernstein-Markov inequality. In this case the zeros of polynomials in $L^{2}(\mu)$ tend to concentrate around the Silov boundary of $K$. In the following we consider the equidistribution of the zeros of polynomials with respect to the Poincar\'e metric at infinity on $\C$.

Of course $\mathbb C$ is a special case of a quasi-projective variety, its complement in $\Proj^{1}$ is the hyperplane at infinity $H_{\infty} = \{z_{0 }=0\}$, via the embedding $\mathbb C \ni \zeta \mapsto [1:\zeta] \in \mathbb P^1$. We denote as usual $U_{j}=\{[z]\in\Proj^{1}:z_{j}\neq0\}$. The hyperplane section bundle $\mO(1)\to\Proj^{1}$ comes along with the canonical (defining) section $s_{0}$ locally given as
$s_{0}|_{U_{0}}=1\cdot e_{0},$ $s_{0}|_{U_{i}} = z_{0}/z_{i} \cdot e_{i}= \xi^0_{i}\cdot e_{i}$, where $e_{j}$ are the canonical frames of $\mO(1)|_{U_{j}}$, $j=0,1$. Now $\cup \{\xi^0_{i}<R\}$ is an open neighbourhood of $H_{\infty}$.\par
Let us consider the charts:
\[\mathbb P^1 = U_0 \cup \{\infty\}=U_{0}\cup U_{1},\quad z:U_0 \to \mathbb C,\quad w:U_1 \to \mathbb C.\]
Consider the divisor $D$ given by the following cover together with meromorphic functions $\{(U_0,1),(U_1,w)\}$. The associated line bundle $[D]$ is defined by the cocycle $U_0\cap U_1, g_{01}=1/w$ and we have $[D]=\mO(1)$. A metric on $[D]$ corresponds to functions $h_i \in \cC^{\infty}(U_i,\R_{>0})$ that satisfy $h_1=|g_{01}|^2h_0.$
In $\{|w|<R\} $ set $h_1=1$ and extend it to a smooth metric over $\mathbb P^1.$ Then $s_{0}=0$ precisely at $\infty$. To determine the $L^2$-condition in the Poincar\'{e} metric it suffices to investigate the integrals in a neighbourhood of $\infty$.
\par
It is well known that the holomorphic sections of $[D]^{\otimes N}$ are identified with complex polynomials of degree $\leqslant N$ in the chart $U_{0}$. We denote this space by $\mathcal H_N$. The Poincar\'{e} metric on $\mathbb C$ is
\[\Theta_{\varepsilon}=\omega_{\FS} -\imat\varepsilon\partial\overline{\partial}\log(-\log\|s_{0}\|^2)^2,\]
the metric on $\mathcal O(1)$ is
\[ h^{\mathcal O(1)}_{\delta}=h^{\mathcal O(1)}\cdot (-\log\|s_{0}\|^2)^{2\delta}.\]
Note that choosing $\delta = 2\pi\varepsilon$ provides $\imat R^{(\mathcal O(1),h^{\mathcal O(1)}_{\delta})}=2\pi\Theta_{\varepsilon}$.

\noindent
A polynomial $P \in \mathcal H_N$ lies in $L^2(\mathbb C,\Theta_{\varepsilon},h^N_{\delta})$ if and only if the integral
\[\int_{|z|>R}{|P(z)|^2\underbrace{(1+|z|^2)^{-N}}_{=h^{\mathcal O(N)}}\big(-\log \left|\tfrac{1}{z}\right|^2\big)^{2N\delta}\underbrace{ \Big\{ \frac{\imat}{2\pi}\frac{dz\wedge d\bar{z}}{(1+|z|^2)^2}-\imat\varepsilon\partial\overline{\partial}\log\big(-\log\left|\tfrac{1}{z}\right|^2\big)^2 \Big\} }_{=\Theta_{\varepsilon}}}\] is finite.
If $\deg\, P=d$, then $|P(z)|^2(1+|z|^2)^{-N}={\rm O}((1+|z|^2)^{-n}),$ with $-2N\leqslant n=2d -2N \leqslant 0$ and in particular bounded.

\noindent
First consider
\[I=\int_{|z|>R}(\log \left|{z}\right|^2)^{2N\delta}\frac{\imat}{2}\frac{dz\wedge d\bar{z}}{(1+|z|^2)^2}\;\cdot\] In polar coordinates
\begin{align*}
I&= C \int_R^{\infty}(\log r^2)^{2N\delta}\frac{rdr}{(1+r^2)^2}\qquad \text{substitute $r^2=e^x$, $2rdr=e^xdx$}\\
&=C^{\prime}\int_{2\log R}^{\infty} x^{2N\delta}\frac{e^xdx}{(1+e^x)^2}\\
&\leqslant C^{\prime}\int\frac{x^{2N\delta}dx}{(1+e^x)}\\
&\leqslant C^{\prime}\int x^{2N\delta}e^{-x}dx <\infty\,.
\end{align*}
We estimate the second integral. We have
\[ \partial\overline{\partial}\log\big(-\log\left|\tfrac{1}{z}\right|^2\big)^2 =\partial\overline{\partial}\log(\log|z|^2)^2 =  \frac{-2dz\wedge d\bar{z}}{|z|^2(\log|z|^2)^2}\,\cdot\]
Therefore
\begin{align*}
\int_{|z|>R}(1+|z|^2)^{-n}(\log\left|{z}\right|^2)^{2N\delta -2} \frac{\imat dz\wedge d\bar{z}}{2|z|^2}&= C \int_R^{\infty} (1+r^2)^{-n}(\log r^2)^{2N\delta -2}\frac{rdr}{r^2}\\
&=C^{\prime}\int_{2\log R}^{\infty}(1+e^x)^{-n}x^{2N\delta-2}dx
\end{align*}
is finite ($N\to\infty$), only if $n< 0.$
This shows that
\[
\mathcal H^N \cap L^2(\mathbb C, \Theta_{\varepsilon}, h^N_{\delta}) = \mathcal H^{N-1}
\]
as sets.
\begin{corollary}
Denote by $\Proj\mathcal H^{N-1}$ the projective space associated to $\mathcal H^{N-1}$. Then
zero-divisors of generic random sequences $(s_{N}) \in \prod_{N=1}^\infty \Proj\mathcal H^{N-1}$
are equidistributed with respect to $\Theta_{\varepsilon}$\,.
\end{corollary}


\def\cprime{$'$} \def\cprime{$'$}
\providecommand{\bysame}{\leavevmode\hbox to3em{\hrulefill}\thinspace}
\providecommand{\MR}{\relax\ifhmode\unskip\space\fi MR }
\providecommand{\MRhref}[2]{%
  \href{http://www.ams.org/mathscinet-getitem?mr=#1}{#2}
}
\providecommand{\href}[2]{#2}

\end{document}